\documentclass[pra,twocolumn]{revtex4}

\usepackage{amssymb,amsmath,amsfonts,latexsym}
\usepackage{dsfont,graphicx,graphics,epsfig,color}
\usepackage{wasysym}
\newtheorem{thm}{Theorem}
\newtheorem{prop}[thm]{Proposition}
\newtheorem{fact}[thm]{Fact}

\newtheorem{cor}[thm]{Corollary}

\def\a{\mathbf{a}}
\def\b{\mathbf{b}}
\def\r{\mathbf{r}}

\begin{document}

\title{Rank penalized estimation of a quantum system}
\author{Pierre Alquier$^{1}$, Cristina Butucea$^{2,4}$, Mohamed Hebiri$^{2}$,
Katia Meziani$^{3}$,
Tomoyuki Morimae$^{5,6}$}
\address{$^1$ University College Dublin, School of Mathematical Sciences,
528 James Joyce Library, Belfield, Dublin 4, Ireland\\
$^2$ LAMA, Universit\'e Marne-la-Vall\'ee, 5, bd. Descartes, Champs-sur-Marne,
77454 Marne-la-Vall\'ee Cedex 2\\
$^3$ CEREMADE, Universit\'e Paris Dauphine, Pl. du Mar\'echal De Lattre De
Tassigny, 75775 Paris Cedex 16 \\
$^4$ CREST, ENSAE, Timbre J-350, 3, av. Pierre Larousse, 92240 Malakoff\\
$^5$ Department of Physics, Imperial College London,
SW7 2AZ, United Kingdom\\
$^6$ ASRLD Unit, Gunma University, 1-5-1 Tenjin-cho Kiryu-shi Gunma-ken, 376-0052, Japan
}
\begin{abstract}
We introduce a new method to reconstruct the density matrix $\rho$ of a
system of $n$-qubits and estimate its rank $d$ from data obtained by quantum state tomography measurements repeated $m$ times.
The procedure consists in minimizing the risk of a linear  estimator
$\hat{\rho}$ of $\rho$ penalized by given rank
(from 1 to $2^n$), where  $\hat{\rho}$ is previously obtained by the moment method. We obtain simultaneously an estimator
of the rank and the resulting density matrix associated to this rank. We establish an upper bound for the error of penalized estimator,
evaluated with the Frobenius norm, which is of order $dn(4/3)^n /m$ and consistency for the estimator of the rank. The proposed methodology is  computationaly efficient and is illustrated with some example states and real experimental
data sets.


\end{abstract}

\maketitle


\section{Introduction}
The experimental study of quantum mechanical systems has made huge progress
recently motivated by quantum information science.
Producing and manipulating many-body quantum mechanical systems have been
relatively easier over the last decade.
One of the most essential goals in such experiments
is to reconstruct quantum states
via quantum state tomography (QST).
The QST is an experimental
process where the system is repeatedly measured with different elements of a
positive operator valued measure (POVM).

Most popular methods for estimating the state from such data are: linear
inversion \cite{VR}, \cite{RMH},
maximum likelihood \cite{BAP}, \cite{JKM}, \cite{L}, \cite{BK2}, \cite{Y}
and Bayesian inference
\cite{AS}, \cite{BBG}, \cite{BK1} (we also refer the reader to \cite{AF+,BD+}
and references therein). Recently,
different approaches brought up-to-date statistical techniques in this field.
The estimators are obtained via minimization of a penalized risk. The
penalization will subject the estimator to constraints. In \cite{Kolt} the
penalty is the Von Neumann entropy of the state, while \cite{GLF}, \cite{Gross}
use the $\mathbb{L}_1$ penalty, also known as the Lasso matrix estimator, under
the assumption that the state to be estimated has low rank. These last papers
assume that the number of measurements must be minimized in order to recover all
the information that we need. The ideas of matrix completion is indeed, that,
under the assumptions that the actual number of underlying parameters is small
(which is the case under the low-rank assumption) only a fraction of all
possible measurements will be sufficient to recover these parameters. The choice
of the measurements is randomized and, under additional assumptions, the
procedure will recover the underlying density matrix as well as with the full
amount of measurements (the rates are within $log$ factors slower than the rates
when all measurements are performed).

\bigskip

In this paper, we suppose that a reasonable amount $m$ (e.g. $m=100$) of data is
available from all possible measurements.
We implement a method to recover the
whole density matrix and estimate its rank from this huge amount of data.
This problem was already considered by Gu\c t\u a, Kypraios and Dryden~\cite{GKD} who propose a maximum likelihood estimator of the state.
Our method is relatively easy to implement and computationally efficient. Its
starting point is a linear estimator obtained by the moment method (also known as the inversion method), which is
projected on the set of matrices with fixed, known rank. A data-driven procedure
will help us select the optimal rank and minimize the estimators risk in
Frobenius norm. We proceed by minimizing the risk of the linear estimator,
penalized by the rank.
When estimating the density matrix of a $n$-qubits system, our final procedure
has the risk (squared Frobenius norm) bounded by  $dn(4/3)^n/m$, where $d$
between 1 and $2^n$ is the rank of the matrix.

The inversion method is known to be computationally easy but less convenient
than constrained maximum
likelihood estimators as it does not produce a density matrix as an output. We
revisit the moment
method in our setup and argue that we can still transform the output into a
density matrix, with
the result that the distance to the true state can only be decreased in the
proper norm.

We shall indicate how to transform the linear estimator into a physical state  with fixed, known rank.
Finally, we shall select the estimator which fits best to the data in terms of a
rank-penalized error. Additionally, the rank selected by this procedure is a
consistent estimator of the true rank $d$ of the density matrix.

We shall apply our procedure to the real data issued from experiments on systems
of
4 to 8 ions.
Trapped ion qubits are a promising candidate for building a quantum computer.
An ion with a single electron in the valence shell is used. Two qubit states are
encoded in two energy levels of the valence electrons, see \cite{BSG}, \cite{GKD}, \cite{Monz14}.

\bigskip

The structure of the paper is as follows. Section 2 gives notation and setup of
the problem. In Section 3 we present the moment method. We first change
coordinates of the density matrix in the basis of Pauli matrices and vectorize the
new matrix. We give properties of the linear operator which takes this vector of
coefficients to the vector of probabilities $p(\a,\r)$. These are the
probabilities to get a certain outcome $\r$ from a given measurement indexed by
$\a$ and that we actually estimate from data at our disposal. We prove the
invertibility of the operator, i.e. identifiability of the model (the
information we measure enables us to uniquely determine the underlying
parameters). Section 4 is dedicated to the estimation procedure. The linear
estimator will be obtained by inversion of the vector of estimated coefficients.
We describe the rank-penalized estimator and study its error bounds. We study
the numerical properties of our procedure on example states and apply them to
experimental real-data in Section 5. The last section is dedicated to proofs.

\section{Basic notation and setup}
\label{section_notations}

We have a system of $n$ qubits. This system is represented by a
$2^{n}\times 2^{n}$ density matrix $\rho$, with coefficients in $\mathbb{C}$.
This matrix is Hermitian $\rho^\dagger=\rho$, semidefinite positive
$\rho\geq 0$
and has ${\rm Tr}(\rho)=1$.
The objective is to estimate $\rho$, from measurements of many independent
systems,
identically prepared in this state.

For each system, the experiment provides random data from separate measurements
of Pauli matrices $\sigma_x, \, \sigma_y, \, \sigma_z$ on each particle.
The collection of measurements which
are performed writes
\begin{equation} \label{mes}
\{\sigma_{\a} = \sigma_{{a}_1} \otimes \ldots \otimes \sigma_{{a}_n}, \quad
\mathbf{a} \in \mathcal{E}^n = \{x,y,z\}^{n}\},
\end{equation}
where $\mathbf{a} = (a_1,\ldots,a_n)$ is a vector taking values in
$\mathcal{E}^n $ which identifies the experiment.

The outcome of the experiment will be a vector
$\r \in\mathcal{R}^n = \{-1,1\}^{n}$. It follows from the basic principles
of quantum mechanics
that the outcome of any experiment indexed by $\mathbf{a}$ is actually a random
variable,
say $R^{\mathbf{a}}$,
and that its distribution is given by:
\begin{equation}
\label{loi}
\forall \mathbf{r} \in \mathcal{R}^n,
   \mathbb{P}(R^\mathbf{a} = \mathbf{r}) = {\rm Tr}
     \left(\rho \cdot P_{r_{1}}^{a_{1}} \otimes \dots \otimes
                    P_{_{r_n}}^{a_{n}} \right),
\end{equation}
where the matrices $P_{r_i}^{a_i}$ denote the projectors on the eigenvectors of
$\sigma_{a_i}$
associated to the eigenvalue $r_i$, for all $i$ from 1 to $n$.

For the sake of simplicity, we introduce the notation
$$ P_{\mathbf{r}}^{\mathbf{a}} := P_{r_{1}}^{a_{1}} \otimes \dots \otimes
                    P_{_{r_n}}^{a_{n}}.$$
As a consequence we have the shorter writing for~\eqref{loi}:
$\mathbb{P}(R^\mathbf{a} = \mathbf{r}) = {\rm Tr} \left(\rho \cdot
 P_{\mathbf{r}}^{\mathbf{a}} \right)$.

The tomographic inversion method for reconstructing $\rho$ is based on
estimating
probabilities $p(\a,\r) := \mathbb{P}(R^\mathbf{a} = \mathbf{r})$ by $\hat
p(\a,\r)$ from available data and solving the linear system of equations
\begin{equation} \label{loihat}
\hat p(\a,\r) = {\rm Tr} \left(\hat {\rho} \cdot
 P_{\mathbf{r}}^{\mathbf{a}} \right).
\end{equation}
It is known in statistics as the method of moments.

\bigskip

We shall use in the sequel the following notation:
$\|A\|_F^2 = {\rm Tr}(A^\dagger
A)$ denotes the Frobenius norm
and $\|A\| = \sup_{v \in \mathbb{R}^{d}, |v|_2 = 1} |A v|_2 $ the operator
sup-norm for any $d \times d$ Hermitian matrix $A$, $|v|_2$ is the Euclidean
norm of the vector $v \in \mathbb{R}^d$.

\bigskip

In this paper, we give an explicit inversion formula for solving (\ref{loi}).
Then, we apply the inversion procedure to equation (\ref{loihat}) and this will
provide us an unbiased estimator $\hat{\rho}$ of $\rho$.
Finally, we project this estimator on the subspace of matrices of rank $k$ ($k$
between 1 and $2^n$)
and thus choose, without any a priori assumption, the estimator which best fits
the data. This is done by minimizing the penalized risk
$$
\|R-\hat{\rho}\|_F^2+\nu \cdot {\rm rank}(R),
$$
where the minimum is taken over all Hermitian, positive semidefinite matrices
$R$. Note that the output is not a proper density matrix. Our last step will transform the output in a physical state.
The previous optimization program has an explicit and easy to implement solution.
The procedure will also estimate the rank of the matrix which best fits data.
We actually follow here the rank-penalized estimation method proposed in the
slightly different problems of matrix regression. This problem recently received
a lot of attention in the statistical community \cite{BSW,Klopp,NW,RT}
and Chapter 9 in \cite{KoltStF}. Here, we follow the computation in \cite{BSW}.

In order to give such explicit inversion formula we first change the coordinates
of the matrix $\rho$ into a vector $\vec{\rho} \in \mathbb{R}^{4n}$ on a
convenient basis.
The linear inversion also gives information about the quality of each estimator
of the coordinates in $\vec{\rho}$. Thus we shall see that we have to perform all
measurements $\sigma_\a$ in order to recover (some) information on each
coordinate of $\vec{\rho}$. Also, some coordinates are estimated from several
measurements and the accuracy of their estimators is thus better.

To our knowledge, this is the first time that rank penalized estimation of a
quantum state is performed.
Parallel work of Gu\c t\u a {\it et al.}~\cite{GKD} addresses the same issue via
the maximum likelihood procedure.
Other adaptive methods include matrix completion for low-rank matrices
\cite{Candes,GLF,Gross,KLT} and for matrices with small Von Neumann entropy
\cite{Kolt}.

\section{Identifiability of the model}

Note the problem of state tomography with mutually unbiased bases,
described in Section~\ref{section_notations}, was considered in
Refs.~\cite{Filippov,Ivonovic}. In this section, we introduce some
notation used throughout the paper, and remind some facts that were
proved for example in~\cite{Ivonovic} about the identifiability of
the model.

A model is identifiable if, for different values of the underlying parameters,
we get different likelihoods (probability distributions) of our sample data.
This is a crucial property for establishing the most elementary convergence
properties of any estimator.

The first step to explicit inversion formula is
to express $\rho$ in the $n$-qubit Pauli basis.
In other words, let us put
$ \mathcal{M}^n = \{I,\, x,\, y,\, z\}^{n}$ and $\sigma_I = I$.  For all $\b \in
\mathcal{M}^n$, denote similarly to (\ref{mes})
\begin{equation}\label{basis}
\{ \sigma_{\b} = \sigma_{b_1} \otimes \ldots \otimes \sigma_{b_n}, \quad \b \in
\mathcal{M}^n \}.
\end{equation}
Then, we have the following decomposition:
$$
\rho = \sum_{\b \in \mathcal{M}^n} \rho_{\b} \cdot \sigma_\b, \quad
\mbox{ with } \rho_{(I,\dots,I)}= \frac 1{2^n}.
$$

We can plug this last equation into~\eqref{loi} to obtain, for
$\mathbf{a}\in\mathcal{E}^{n}$
and $\mathbf{r}\in\mathcal{R}^n$,
\begin{align*}
& \mathbb{P}(R^\mathbf{a} = \mathbf{r})  = {\rm Tr} \left(\rho \cdot
 P_{\mathbf{r}}^{\mathbf{a}} \right)
\\
& = {\rm Tr} \left( \sum_{\b \in \mathcal{M}^n}
\rho_{\b} \cdot \sigma_{\b} \cdot P_{\mathbf{r}}^{\mathbf{a}} \right)
\\
& = \sum_{\b\in\mathcal{M}^n}
\rho_{\b}  {\rm Tr} \left( \left(\sigma_{b_1} \otimes \dots \otimes
\sigma_{b_n}\right)
  \left(P_{r_{1}}^{a_{1}} \otimes \dots \otimes
                    P_{_{r_n}}^{a_{n}} \right) \right)
\\
& = \sum_{\b\in\mathcal{M}^n}
\rho_{\b} \prod_{j=1}^{n} {\rm Tr}(\sigma_{b_j} P_{_{r_j}}^{a_{j}}).
\end{align*}
Finally, elementary computations lead to ${\rm Tr}(I P_{s}^{t})=1$ for
any $s\in\{-1,1\}$ and $t\in\{x,y,z\}$, while ${\rm Tr}(\sigma_t P_{s}^{t'})
= s \delta_{t,t'}$ for any $s\in\{-1,1\}$, $(t,t')\in\{x,y,z\}^2$ and $\delta$ denotes the Kronecker symbol.

For any $\mathbf{b}\in\mathcal{M}^{n}$, we denote by $E_\b = \{j \in
\{1,\ldots,n\}: b_j=I\}$.
The above calculation leads to the following fact,
which we will use later.

\begin{fact}
\label{propprob}
For $\mathbf{a}\in\mathcal{E}^{n}$,
and $\mathbf{r}\in\mathcal{R}^n$, we have
$$ \mathbb{P}(R^\mathbf{a} = \mathbf{r})
= \sum_{\b \in \mathcal{M}^n}
\rho_{\mathbf{b}} \cdot \prod_{j \not\in E_\b} r_j \, \mathbf{I}(a_j=b_j).
$$
\end{fact}
Let us consider, for example, $\b = (x,\ldots,x)$, then the associated set $E_{\b}$ is empty
and $\mathbb{P}(R^{(x,\ldots,x)}=\r)$ is the only probability depending
on $\rho_{(x,\ldots,x)} $ among other coefficients. Therefore, only the
measurement $(\sigma_x,\ldots,\sigma_x)$ will bring information on this
coefficient. Whereas, if $\b = (I,I,x,\ldots,x)$, the set $E_\b$ contains 2
points. There are $3^2$ measurements $\{ (\sigma_x, ...,\sigma_x)$, ...,
$(\sigma_z,\sigma_z, \sigma_x,...,\sigma_x) \}$ that will bring partial
information on $\rho_\b$. This means, that a coefficient $\rho_{\b}$ is
estimated with higher accuracy as the size of the set $E_\b$ increases.

For the sake of shortness, let us put in vector form:
\begin{eqnarray*}
\vec{\rho}  &:= & (\rho_{\b})_{\b \in\mathcal{M}^n}\\
\mathbf{p} & :=&
\left(p_{(\mathbf{r},\mathbf{a})}\right)_{(\mathbf{r},\mathbf{a})\in(\mathcal{R}
^n\times\mathcal{E}^{n})}
 = (\mathbb{P}(R^\mathbf{a} = \mathbf{r}))
_{(\mathbf{r},\mathbf{a})\in(\mathcal{R}^n\times\mathcal{E}^{n})}.
\end{eqnarray*}
Our objective is to study the invertibility of the operator
\begin{align*}
          \mathbb{R}^{4^n} & \rightarrow
                      \mathbb{R}^{6^n} \\
          \vec{\rho}
& \mapsto
\mathbf{p}.
\end{align*}

Thanks to Fact~\ref{propprob}, this operator is linear. It can then be
represented by a matrix
$\mathbf{P} = [P_{(\mathbf{r},\mathbf{a}),\b}]_{
(\mathbf{r},\mathbf{a})\in(\mathcal{R}^n\times \mathcal{E}^{n}),\b
\in\mathcal{M}^n}$,
we will then have:
\begin{equation}
\label{operator}
  \forall (\mathbf{r},\mathbf{a})\in (\mathcal{R}^n\times \mathcal{E}^{n}),
\quad p_{(\mathbf{r},\mathbf{a})} = \sum_{\b\in\mathcal{M}^n}
         \rho_{\b} P_{(\mathbf{r},\mathbf{a}),\b}
\end{equation}
and from Fact~\ref{propprob} we know that
$$
P_{(\mathbf{r},\mathbf{a}),\b} =
         \prod_{j\not\in E_\b} r_j \, \mathbf{I}{(a_j=b_j)} .
$$
We want to solve the linear equation $\mathbf{P} \vec{\rho} = \mathbf{p}$.
Recall that $E_\b$ is the set of indices where the vector $\b$ has an $I$
operator. Denote by $d(\b)$ the cardinality of the set $E_\b$.
\begin{prop}
\label{inversion}
The matrix $ \mathbf{P}^{T} \mathbf{P} $ is a diagonal matrix with non-zero
coefficients given by
$$
(\mathbf{P}^{T} \mathbf{P})_{\b,\b} = 3^{d(\b)} \, 2^n.
$$
As a consequence the operator is invertible, and the equation $\mathbf{P} \vec{\rho} =
\mathbf{p}$ has a unique solution:
$$ \vec{\rho} = (\mathbf{P}^{T} \mathbf{P})^{-1} \mathbf{P}^{T} \mathbf{p} .$$
\end{prop}
In other words, we can reconstruct $\vec{\rho}=(\rho_{\b})_{\b\in\mathcal{M}^{n}}$
from $\mathbf{p}$,
in the following way:
$$ \rho_{\b} = \frac{1}{3^{d(\b)} 2^n}
\sum_{(\mathbf{r}, \a)\in (\mathcal{R}^n\times \mathcal{E}^{n}) }
p_{(\mathbf{r}, \mathbf{a})}
\mathbf{P}_{(\mathbf{r},\mathbf{a}),\b}.$$
This formula confirms the intuition that, the larger is $d(\b)$, the more
measurements $\sigma_\a$ will contribute to recover the coefficient $\rho_\b$.
We expect higher accuracy for estimating $\rho_\b$ when $d(\b)$ is large.

\section{Estimation procedure and error bounds}

In practice, we do not observe $\mathbb{P}(R^\mathbf{a} = \mathbf{r})$
for any $\mathbf{a}$ and $\mathbf{r}$.
For any $\mathbf{a}$, we have a set of $m$ independent experiments, whose
outcomes are
denoted by $R^{\a,i}$, $1\leq i \leq m$. Our setup is that the
$R^{\mathbf{a},i}$
are independent, identically distributed (i.i.d.) random variables, distributed
as $R^\mathbf{a}$.

We then have a natural estimator for $p_{(\mathbf{r},\mathbf{a})} =
\mathbb{P}(R^\mathbf{a} = \mathbf{r})$:
$$ \hat{p}_{(\mathbf{r},\mathbf{a})} =
   \frac{1}{m}\sum_{i=1}^{m} \delta_{R^{\mathbf{a},i},\mathbf{r}}. $$
We can of course write $\hat{\mathbf{p}} =
(\hat{p}_{(\mathbf{r},\mathbf{a})})_{(\mathbf{r},\mathbf{a})\in(\mathcal{R}
^n\times \mathcal{E}^{n})}
$.

\subsection{Linear estimator}

We apply the inversion formula to the estimated vector $\hat{\mathbf{p}}$.
Following Proposition~\ref{inversion} we can define:
\begin{equation} \label{vechat}
\hat{\vec{\rho}} = (\mathbf{P}^{T} \mathbf{P})^{-1} \mathbf{P}^{T}
\hat{\mathbf{p}}.
\end{equation}
Put it differently:
$$ \hat{\rho}_{\b} = \frac{1}{3^{d(\b)} 2^n}
\sum_{(\mathbf{r},\mathbf{a})\in(\mathcal{R}^n\times \mathcal{E}^{n})}
    \hat{p}_{(\mathbf{r},\mathbf{a})}
 \mathbf{P}_{(\mathbf{r},\mathbf{a}),\b} $$
and then, the linear estimator obtained by inversion, is
\begin{equation} \label{mathat}
\hat{\rho} = \sum_{\b \in \mathcal{M}^n} \hat{\rho}_{\b}
 \sigma_{\b}.
 \end{equation}

The next result gives asymptotic properties of the estimator
$\hat{\vec{\rho}}$ of $\vec{\rho}$.
 \begin{prop}\label{lambda}
 The estimator $\hat{\vec{\rho}}$ of $\vec{\rho}$, defined in (\ref{vechat}) has the
following properties:
\begin{enumerate}
\item it is unbiased, that is $\mathbb{E}[\hat{\vec{\rho}}] = \vec{\rho}$;
\item it has variance bounded as follows
$$
Var(\hat \rho_{\b}) \leq \frac 1{3^{d(\b)} 4^n m};
$$
\item  for any $\varepsilon>0$,
$$
\mathbb{P}\left( \left\| \hat{\rho} - \rho \right\|^{2} \geq
4\sqrt{2 \left(\frac{4}{3}\right)^n   \frac{n\log(2) - \log(\varepsilon) }{m} }\right) \leq \varepsilon.
$$
\end{enumerate}
 \end{prop}
Note again that the accuracy for estimating $\rho_\b$ is higher when $d(\b)$ is
large. Indeed, in this case more measurements bring partial information on $\rho_\b$.

The concentration inequality gives a bound on the norm
$|\hat{\vec{\rho}} - \vec{\rho}|_\infty$ which is valid with high probability.
This quantity is related to $\|\hat{\rho} - \rho\|$ in a way that will be explained later on.
The bound we obtain above depends on $\log(2^n)$,
which is expected as $4^n-1$ is the total number of parameters of a full rank system.
This factor appears in the Hoeffding inequality that we use in order to prove this bound.

\subsection{Rank penalized estimator}

We investigate low-rank estimates of $\rho$ defined in (\ref{mathat}).
From now on, we follow closely the results in \cite{BSW} which were obtained for
a matrix regression model, with some differences as our model is different.
Let us, for a positive real value $\nu$ study the estimator:
\begin{equation}\label{pen}
\hat{\rho}_{\nu} = \arg\min_{R} \left[
\left\|R -  \hat{\rho} \right\|_F^{2} + \nu \cdot {\rm rank}(R)
\right] ,
\end{equation}
where the minimum is taken over all Hermitian matrices $R$.
In order to compute the solution of this optimization program, we may write it in a more convenient form since
\begin{align}\label{penRank}
&\min_{R} \left[
\left\|R -  \hat{\rho} \right\|_F^{2} + \nu \cdot {\rm rank}(R)
\right]\nonumber \\
&= \min_k \min_{R: {\rm rank}(R)=k} \left[
\left\|R -  \hat{\rho} \right\|_F^{2} + \nu \cdot k
\right] .
\end{align}

An efficient algorithm is available to solve the minimization program~\eqref{penRank} as a spectral-based
decomposition algorithm provided in \cite{RV98}. Let us denote by $\hat R_k$ the
matrix such that $\|\hat R_k -  \hat{\rho}\|_F^2 = \min_{R: {\rm rank}(R)=k} \left[
\left\|R -  \hat{\rho} \right\|_F^{2} + \nu \cdot k \right]$. This is a projection of the linear estimator on the space of matrices with fixed (given) rank $k$. Our procedure selects automatically out of data the rank $\hat k$. We see in the sequel that the estimators $\hat R_{\hat k}$ and $\hat{\rho}_\nu$ actually coincide.

We study the
statistical performance from a numerical point of view later on.
\begin{thm}\label{res}
For any $\theta >0$ put $c(\theta) = 1+2/\theta$. We have on the event
$\{\nu \geq (1+\theta) \|\hat{\rho}-\rho\|^2\}$ that
\begin{equation*}\label{penres}
\|\hat{\rho}_{\nu} - \rho\|_F^2
\leq \min _k\left\{
c^2(\theta) \sum_{j>k} \lambda_j^2(\rho) + 2 c(\theta) \nu k
\right\},
\end{equation*}
where $\lambda _j(\rho)$ for $j=1,\ldots ,2^n$ are the eigenvalues of
$\rho$ ordered
decreasingly.
\end{thm}
Note that, if ${\rm rank}(\rho) = d$, for some $d$ between 1 and $2^n$, then the
previous inequality becomes
$$
\|\hat{\rho}_{\nu} - \rho\|_F^2
\leq 2 c(\theta) \nu d
.
$$
%
%
Let us study the choice of $\nu $ in Theorem~\ref{res} such that the
probability of the event
$\{ \nu \geq (1+\theta) \|\hat{\rho}-\rho\|^2 \}$ is small.
By putting together the previous theorem and Proposition~\ref{lambda},
we get the following result:
\begin{cor}
\label{CorBornNuData}
For any $\theta >0$ put $c(\theta) = 1+2/\theta$ and for some small
$\varepsilon >0$ choose
$$
\nu (\theta, \varepsilon)  =
32(1+\theta) \left(\frac{4}{3}\right)^n\frac{n\log(2) - \log(\varepsilon)}{m}
$$
Then, we have
$$
\|\hat{\rho}_{\nu(\theta, \varepsilon)} - \rho\|_F^2
\leq \min _k\left\{
c^2(\theta) \sum_{j>k} \lambda_j^2(\rho) + 2 c(\theta) \nu k
\right\},
$$
with probability larger than $1-\varepsilon$.

\end{cor}
Again, if the true rank of the underlying system is $d$, we can write that, for any $\theta>0$ and for some small $\varepsilon >0$:
$$
\|\hat{\rho}_{\nu} - \rho\|_F^2
\leq 64c(\theta)(1+\theta)d \left(\frac{4}{3}\right)^n\frac{n\log(2) - \log(\varepsilon)}{m} ,
$$
with probability larger than $1-\varepsilon$. If $\|\cdot \|_1$ denotes the trace norm of
a matrix, we have $\|M\|_{1}\leq 2^{\frac{n}{2}} \|M\|_{F}$ for any matrix $M$ of size
$2^n \times 2^n$. So, we deduce from the previous bound that
$$
\|\hat{\rho}_{\nu} - \rho\|_1^2
\leq 64c(\theta)(1+\theta)d \left(\frac{8}{3}\right)^n\frac{n\log(2) - \log(\varepsilon)}{m}.
$$

The next result will state properties of $\hat k$, the rank of the final
estimator $\hat{\rho}_\nu$.
\begin{cor} \label{hatrank} If there exists $k$ such that $\lambda_k (\rho)
>(1+\delta) \sqrt{\nu}$ and $\lambda_{k+1} (\rho) < (1-\delta)
\sqrt{\nu}$ for some
$\delta $ in $(0,1]$, then
$$
\mathbb{P} (\hat k = k) \geq 1- \mathbb{P} (\|\hat{\rho} - \rho\|\geq
\delta \sqrt{\nu}).
$$
\end{cor}
From an asymptotic point of view, this corollary means that, if $d$ is the rank
of the underlying matrix $\rho$, then our procedure is consistent in
finding the rank as the number $m$ of data per measurement increases.
Indeed, as $\sqrt{\nu}$ is an upper bound of the norm
$\|\hat{\rho} - \rho\|$, it tends to 0 asymptotically and therefore
the assumptions of the previous corollary will be checked for $k=d$.
With a finite sample, we deduce from the previous result that $\hat k$ actually
evaluates the first eigenvalue which is above a threshold related to the largest
eigenvalue of the noise $\hat{\rho} - \rho$.

\section{Numerical performance of the procedure}

In this section we implement an efficient procedure to solve the optimization problem~\eqref{penRank} from the previous section.
Indeed, the estimator $\hat \rho$ will be considered as an input from now on. It is computed very efficiently via linear operations and the real issue here is how to project this estimator on a subspace of matrices with smaller unknown rank in an optimal way.
We are interested in two aspects of the method: its ability to select the rank correctly and
the correct choice of the penalty.
First, we explore the penalized procedure on example data and tune the parameter $\nu$ conveniently. In this way, we evaluate the performance of the linear estimator and of the rank selector.
We then apply the method on real data sets.

The algorithm for solving \eqref{penRank} is given in \cite{RV98}. We adapt it to our context and obtain the simple procedure.


\noindent{\bf Algorithm}:\\
\underline{Inputs}: The linear estimator $\hat{\rho}$ and a positive value of the tuning parameter $\nu$\\
\underline{Outputs}: An estimation $\hat k$ of the rank and an approximation $\hat R_{\hat k}$ of the state matrix.
\begin{enumerate}
\item[Step 1.] Compute the eigenvectors $V=[v_1,\ldots,v_{2^n}]$ corresponding to the eigenvalues of the matrix $ \hat{\rho}^{\star}\hat{\rho}$ sorted in decreasing order.
\item[Step 2.] Let $U=\hat{\rho} V$.
\item[Step 3.] For $k=1,\ldots,2^n$, let $V_k$ and $U_k$ be the restrictions to their $k$ first columns
of $V$ and $U$, respectively.
\item[Step 4.] For $k=1,\ldots,2^n$, compute the estimators $\hat R_k = U_k V_k^\star$.
\item[Step 5.] Compute the final solution $\hat R_{\hat k}$, where, for a given positive value $\nu$, $\hat k$ is defined as the minimizer in $k$ over $\{1,\ldots,2^n\}$ of
$$
\left\|\hat R_k -  \hat{\rho} \right\|_F^{2} + \nu \cdot k.
$$
\end{enumerate}

\medskip

The constant $k$ in the above procedure plays the role of the rank and then $\hat R_k$ is the best
approximation of $ \hat{\rho}$ with a matrix of rank $k$.
As a consequence, this approach provides an estimation of both of the matrix
$\rho$ and
of its rank $d$ by $\hat R_{\hat{k}}$ and $\hat{k}$, respectively.

Obviously, this solution is strongly related to the value of the tuning parameter $\nu$. Before dealing
with how to calibrate this parameter, let us present a property that should help us
to reduce the computational cost of the method.

The above algorithm is simple but requires the computation of $2^n$ matrices in Step 3 and Step 4.
We present here an alternative which makes possible to compute only the matrix $\hat R_k $
that corresponds to $k = \hat k$,
and then reduce the storage requirements.

Remember that $\hat k$ is the value of $k$ minimizing the quantity in Step 5 of the above algorithm.
Let $\lambda_1(\hat{\rho})> \lambda_2(\hat{\rho})>...$ be the ordered eigenvalues
of $\sqrt{\hat{\rho}^\star \hat{\rho}}$.
According to \cite[Proposition 1]{BSW}, it turns out that $\hat k$ is the largest $k$ such that
the eigenvalue $\lambda_k(\hat{\rho})$ exceeds the threshold $ \sqrt{\nu}$:
\begin{equation}
\label{RankSelec}
	\hat k = \max \{k : \lambda_k(\hat{\rho}) \geq \sqrt{\nu}\}.
\end{equation}
As a consequence, one can compute the eigenvalues of the matrix
$\sqrt{\hat{\rho}^\star \hat{\rho}}$ and
set $\hat k$ as in~\eqref{RankSelec}.
This value is then used to compute the best solution $\hat R_{\hat k}$ thanks to Step 1 to
Step 4 in the above algorithm, with the major difference that we restrict Step 3 and Step 4
to only $k=\hat k$.

\begin{center}
Example Data
\end{center}

We build artificial
density matrices $\rho$ with a given rank $d$ in $\{1,\ldots,6\}$. These matrices
are $2^n\times 2^n$ with $n=4$ and 5.
To construct such a matrix, we take $\rho$ as $D_d = \frac 1d diag(1 ... 1 0 ... 0)$,  the
diagonal matrix with its first $d$ diagonal terms equal $1/d$, whereas the others
equal zero.

We aim at testing how often we select the right rank based on the method illustrated
in~\eqref{RankSelec} as a function of the rank $d$, and of the number $m$ of repetitions
of the measurements we have in hand.
Our algorithm depends on the tuning parameter $\nu$. We use and compare two different
values of the threshold $\nu$:
denote by $\nu_n^{(1)}$ and $\nu_n^{(2)}$ the values the parameter $\nu$
provided in Theorem~\ref{res} and Corollary~\ref{CorBornNuData} respectively.
That is,
\begin{equation}
\label{TwoParm}
\nu_n^{(1)} =
\Vert\hat{\rho } -  \rho \Vert^2
\quad \text{ and } \quad
\nu_n^{(2)} =   32 (1+\theta) \left(\frac{4}{3}\right)^{n}
\frac{n\log(2)}{m}.
\end{equation}
As established in Theorem~\ref{res}, if the tuning parameter $\nu$ is of order of
the parameter $\nu_n^{(1)} $, the solution of our algorithm is an accurate estimate
of $\rho$. We emphasize the fact that $\nu_n^{(1)} $ is nothing but the estimation error
of our linear estimator $\hat{\rho } $. We study this error below.
On the other hand, the parameter $\nu_n^{(2)}$ is an upper bound of $\nu_n^{(1)}$
that ensures that the accuracy of estimation remains valid with high probability ({\it cf.}  Corollary~\ref{CorBornNuData}).
The main advantage of $\nu_n^{(2)}$ is that it is completely known by the practitioner, which is not the case of $\nu_n^{(1)}$.

\bigskip

\begin{figure}[htp]
\centering
\caption{
\footnotesize{(Color online). Frequency of good selection of the true rank $d$,
based on~\eqref{RankSelec} with $\nu = \nu_n^{(1)} $ (green squares) and
with $\nu = \nu_n^{(2)} $ (blue stars). The results are established on $20$ repetitions.
A value equal to $1$ in the $y$-axis means that the method always selects the good rank, whereas $0$ means that it always fails.
First: $m=50$ measurements -- Second: $m=100$ measurements
}}
\includegraphics[width=0.35 \textwidth]{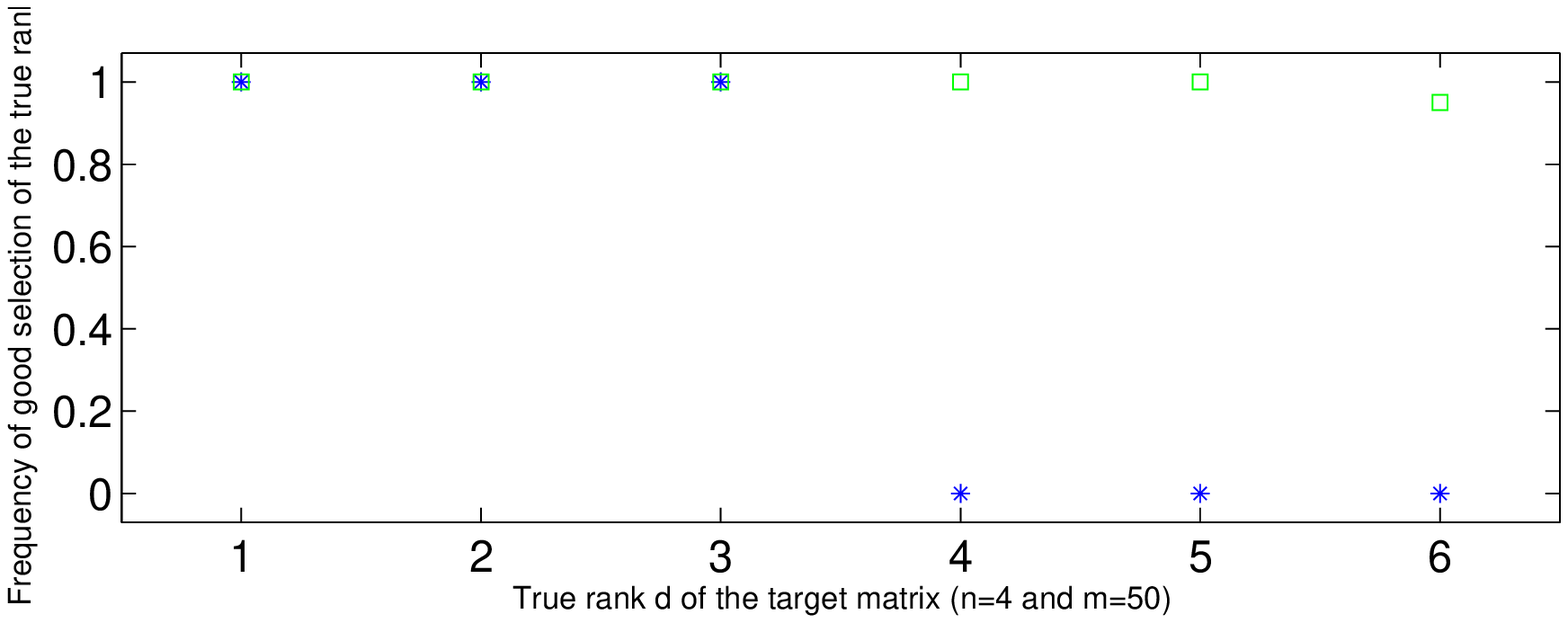}
\includegraphics[width=0.35 \textwidth]{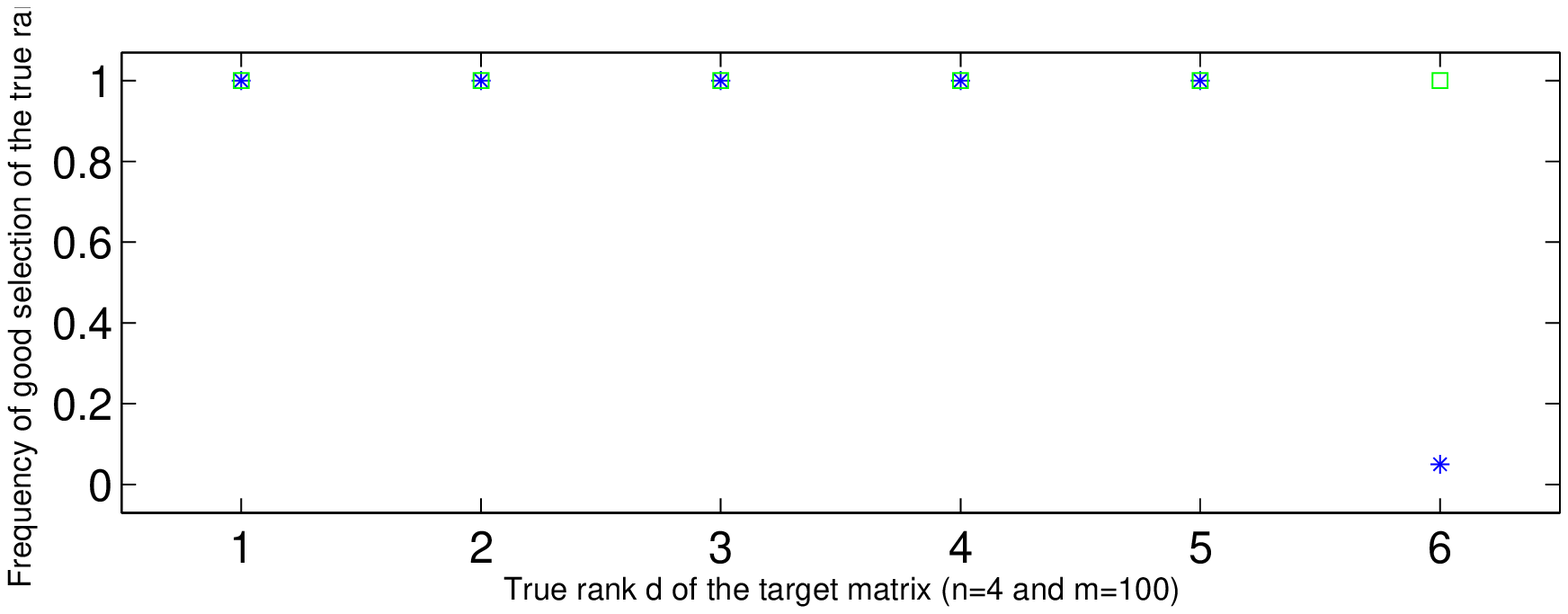}
\label{figErrorVsRank}
\end{figure}

\noindent {\bf Rank estimation.}
Our first goal consists in illustrating the estimation power of our method
in selecting the true rank $d$ based on the calibrations of $\nu$
given by~\eqref{TwoParm}.
We provide some conclusions on the
number of repetitions $m$ of the measurements needed to recover the right rank as a function of this rank.
Figure~\ref{figErrorVsRank} illustrates the evolution of the selection power of our method
based on $\nu_n^{(1)} $ (blue stars) on the one hand, and based on $ \nu_n^{(2)}$ (green squares)
on the other hand.

Two conclusions can be made. First, the method based on
$\nu_n^{(1)} $ is powerful. It almost always selects the right rank. It outperforms
the algorithm based on $\nu_n^{(2)} $. This is an interesting observation. Indeed,
$\nu_n^{(2)} $ is an upper bound of $\nu_n^{(1)} $. It seems that this bound is too large
and can be used only for particular settings. Note however that in the variable selection literature,
the calibration of the tuning parameter is a major issue and is often fixed by Cross-Validation (or other
well-known methods). We have chosen here to illustrate only the result based on our theory
and we will provide later an instruction to properly calibrate the tuning parameter $\nu$.

The second conclusion goes in the direction of this instruction.
As expected, the selection power of the method
(based on both $\nu_n^{(1)} $ and $\nu_n^{(2)} $) increases when the number of
repetition $m$ of the measurements increases. Compare the figure for $m=50$ repetitions to the figure for $m=100$ repetitions in Figure~\ref{figErrorVsRank}. Moreover, for ranks smaller than some values, the
methods always select the good rank. For larger ranks, they perform poorly.
For instance with $m=50$ (a small number of measurements),
we observe that the algorithm based on $\nu_n^{(2)} $ performs poorly
when the rank $d\geq 4$, whereas the algorithm based on $\nu_n^{(1)} $ is still
excellent.\\
Actually, the bad selection when $d$ is large does not mean that the methods perform poorly.
Indeed our definition of the matrix $\rho$ implies that the eigenvalues of the matrix
decrease with $d$. They equal to $1/d$. Therefore, if $\sqrt{\nu}$ is of the same order as
$1/d$, finding the exact rank becomes difficult since this calibration suggests that
the eigenvalues are of the same order of magnitude as the error.
Hence, in such situation, our method adapts to
the context and find the effective rank of $\rho$.
As an example, let consider our study with $n=4$, $m=50$ and $d=6$.
Based on $20$ repetitions of the experiment, we obtain a maximal value of
$\nu_n^{(1)}= \Vert\hat{\rho } -  \rho \Vert^2   $ equal to $0.132$. This value
is quite close to $0.167$, the value of the eigenvalues of $\rho$.
This explains the fact that our method based on $\nu_n^{(1)}$ failed in one iteration
(among $20$) to find the good rank. In this context $\nu_n^{(2)}$ is much larger than
$0.167$ and then our method does not select the correct rank with this calibration in
this setting.
\\
Let us also mention that we explored numerous experiments with other choices
of the density matrix $\rho$. The same conclusion remains valid. When the
error of the linear estimator $\hat{\rho }$ which is given by
$\nu_n^{(1)}= \Vert\hat{\rho } -  \rho \Vert^2   $ is close to the square of the
smallest eigenvalue of $\rho$, finding the exact rank is a difficult task.
However, the method based on $\nu_n^{(1)}$ is still good, but fails sometimes.
We produced data from physically meaningful states: the GHZ-state and the W-state for $n=4$ qubits, as well as a statistical mixture $M_{d,p} = p*GHZ +(1-p)*D_d$, for $d=3$ and $p=0.2$ Note that the rank of $M_{d,p}$ is 4.

\begin{figure}[htp]
\caption{
\footnotesize{(Color online). Evaluation of the operator norm $\sqrt{ \nu_n^{(1)}} =
\Vert\hat{\rho } -  \rho \Vert   $.
The results are established on $20$ repetitions.
Above: $n=4$, $m=50$ repetitions of the measurements ; we compare the errors when $d$ takes values betwenn $1$ and $6$ --
Middle: $n=5$, $m=100$ ; we compare the errors when $d$ takes values between $1$ and $6$ --
Below: the rank equals $d=4$ and compare the error for $m=50$ and $100$.
}}
\includegraphics[width=0.35\textwidth]{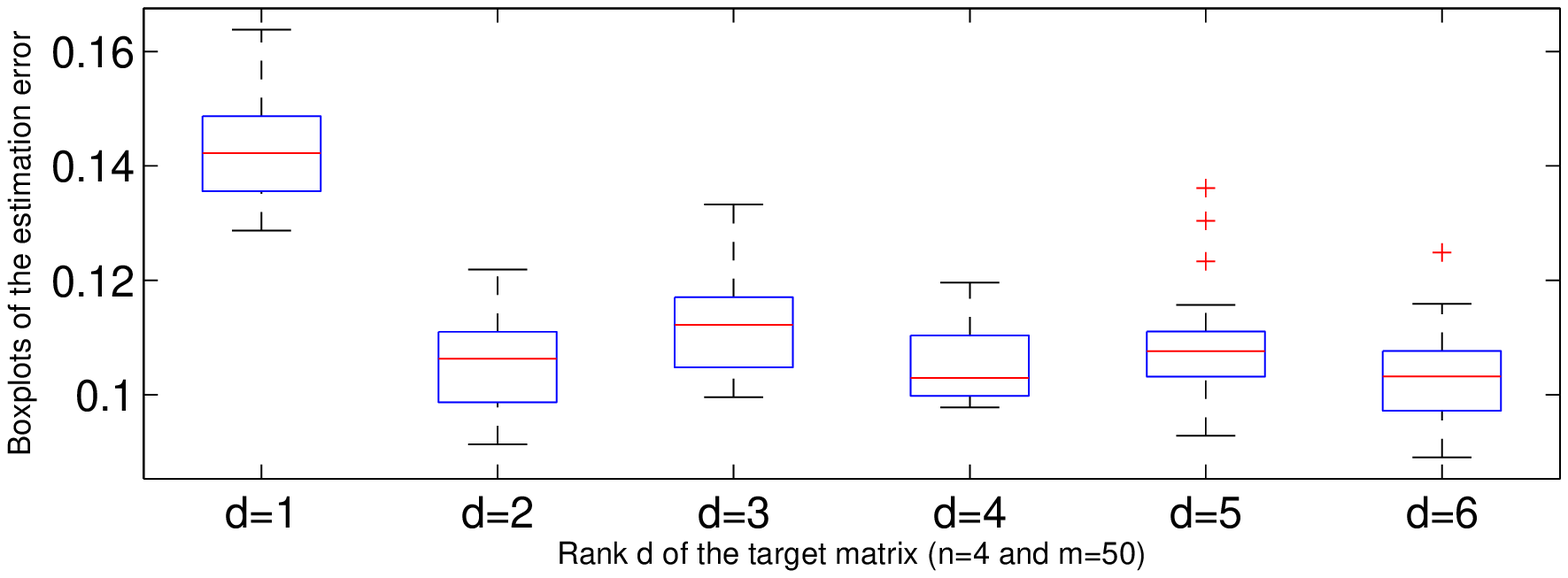}
\includegraphics[width=0.35\textwidth]{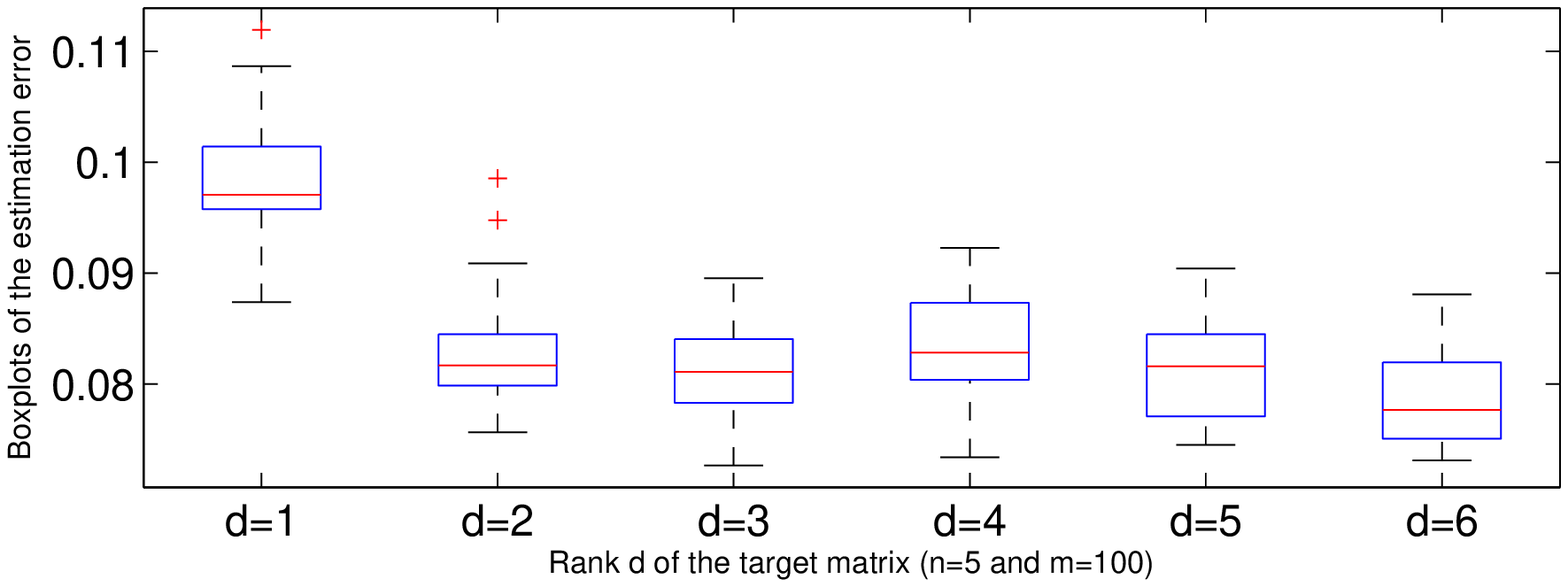}
\includegraphics[width=0.35\textwidth]{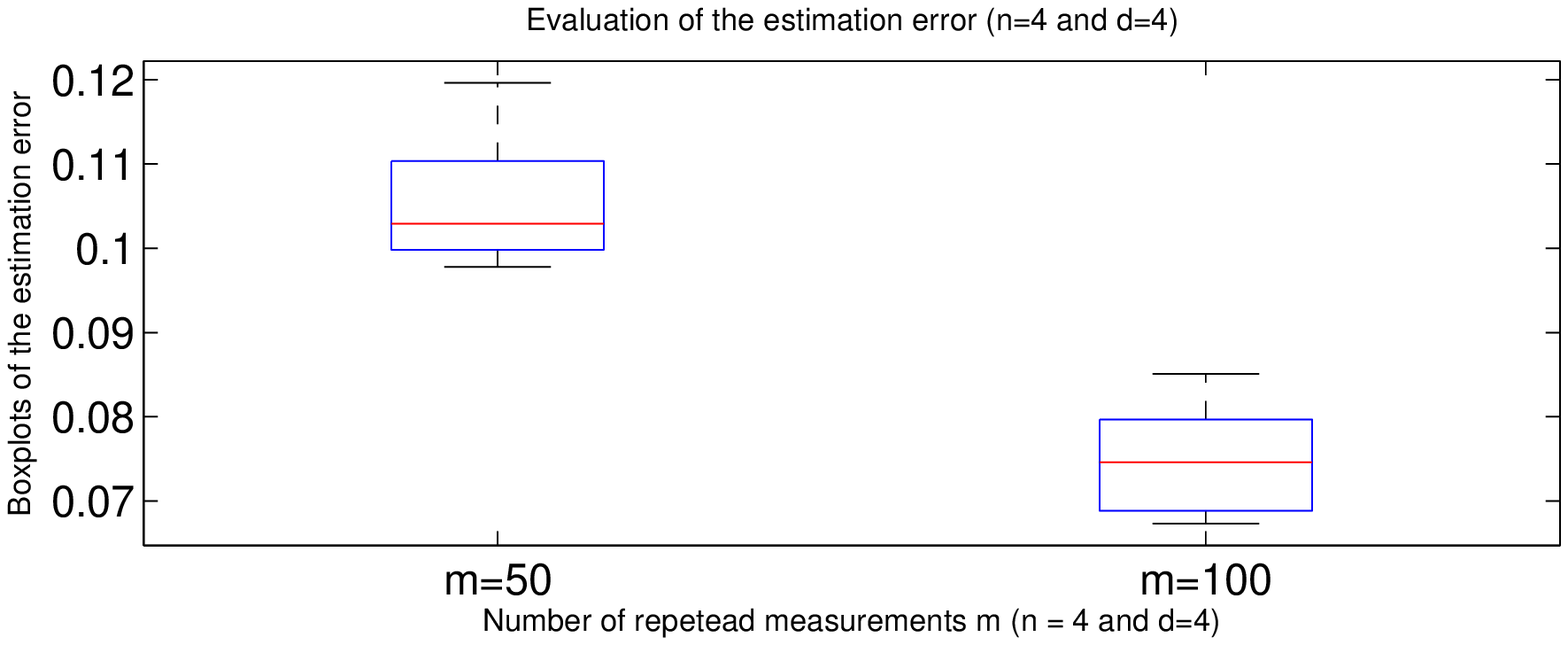}
\label{figOpNorm}
\end{figure}

\bigskip

\noindent {\bf Calibration of the tuning parameter $\nu$.}
The quantity $\nu_n^{(1)}= \Vert\hat{\rho } -  \rho \Vert^2   $ seems to be
very important to provide a good
estimation of the rank $d$ (or more precisely of the effective rank).
Then it is interesting to observe how this quantity behaves.
Figure~\ref{figOpNorm} (Above $m=50$ and $d=4$, and Middle $m=100$ and $d=5$) illustrates how $\nu_n^{(1)} $ varies when the rank increases.
Except for $d=1$, it seems that the value of $\nu_n^{(1)} $ is quite stable. These graphics
are obtained with particular values of the parameters $m$ and $d$, but similar illustrations can be obtained if these parameters change.\\
The main observation according to the parameter $\nu$ is that it decreases with $m$
(see Figure~\ref{figOpNorm} - Below) and is actually independent of the rank $d$ (with
some strange behavior when $d=1$). This is in accordance with the definition
of $\nu_n^{(2)} $ which is an upper bound of $\nu_n^{(1)} $.

\begin{center}
Real-data analysis
\end{center}

In the next paragraph, we propose a 2-steps instruction for practitioners to use our method
in order to estimate a matrix $\rho$ (and its rank $d$) obtained from the data $R^{\a,i}$ we have in
hand with $\a\in\{x,y,z\}$ and $ i\in \{1,\ldots,m\}$.

\bigskip

\noindent {\bf Real Data Algorithm:}\\
\underline{Inputs}: for any measurement $\a\in\{x,y,z\}$ we observe $R^{\a,i}, \ i=1,\ldots,m$.\\
\underline{Outputs}: $\hat k$ and $\hat{R}_{\hat k}$, estimations of the rank $d$ and $\rho$ respectively.
\\
The procedure starts with the linear estimator $\hat{\rho}$ and consists in two steps:

\vspace{2mm}

\noindent  {\it Step A.} Use  $\hat{\rho}$ to simulate repeatedly data with the same parameters $n$ and $m$ as the original problem. Use the data to compute synthetic linear estimators and the mean operator norm of these estimators. They provide an evaluation of the tuning parameter $\tilde\nu_n^{(1)}$.

\vspace{2mm}

\noindent  {\it Step B.} Find $\hat k$ using \eqref{RankSelec} and construct $\hat{R}_{\hat k}$.

\bigskip

\begin{figure}[htp]
\caption{
\footnotesize{(Color online). Eigenvalues of the linear estimator in
increasing order and the penalty choice;
$m=100$ and $n=4, \, 5 $ or 6, respectively.
}}
\includegraphics[width=0.35\textwidth]{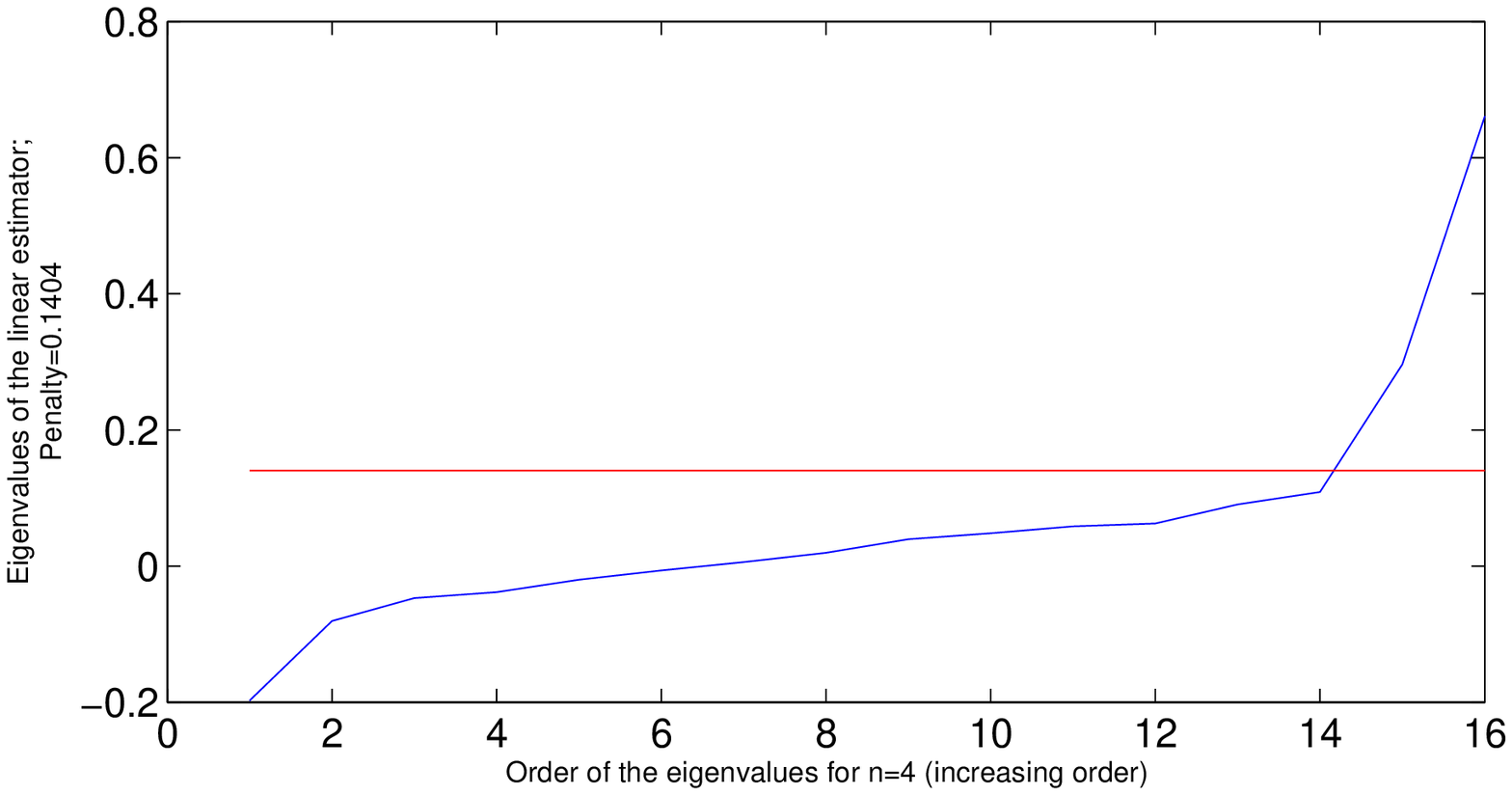}
\includegraphics[width=0.35\textwidth]{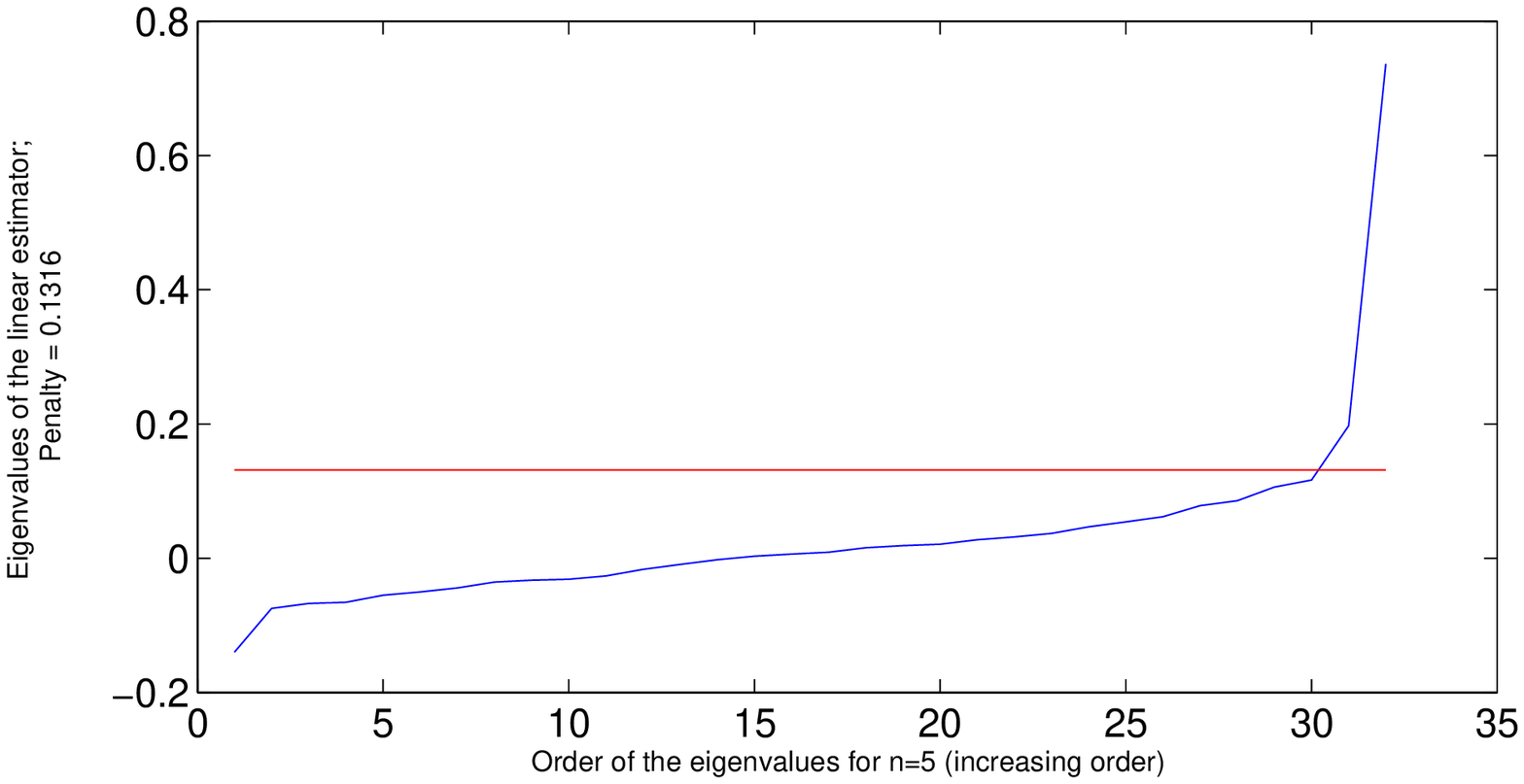}
\includegraphics[width=0.35\textwidth]{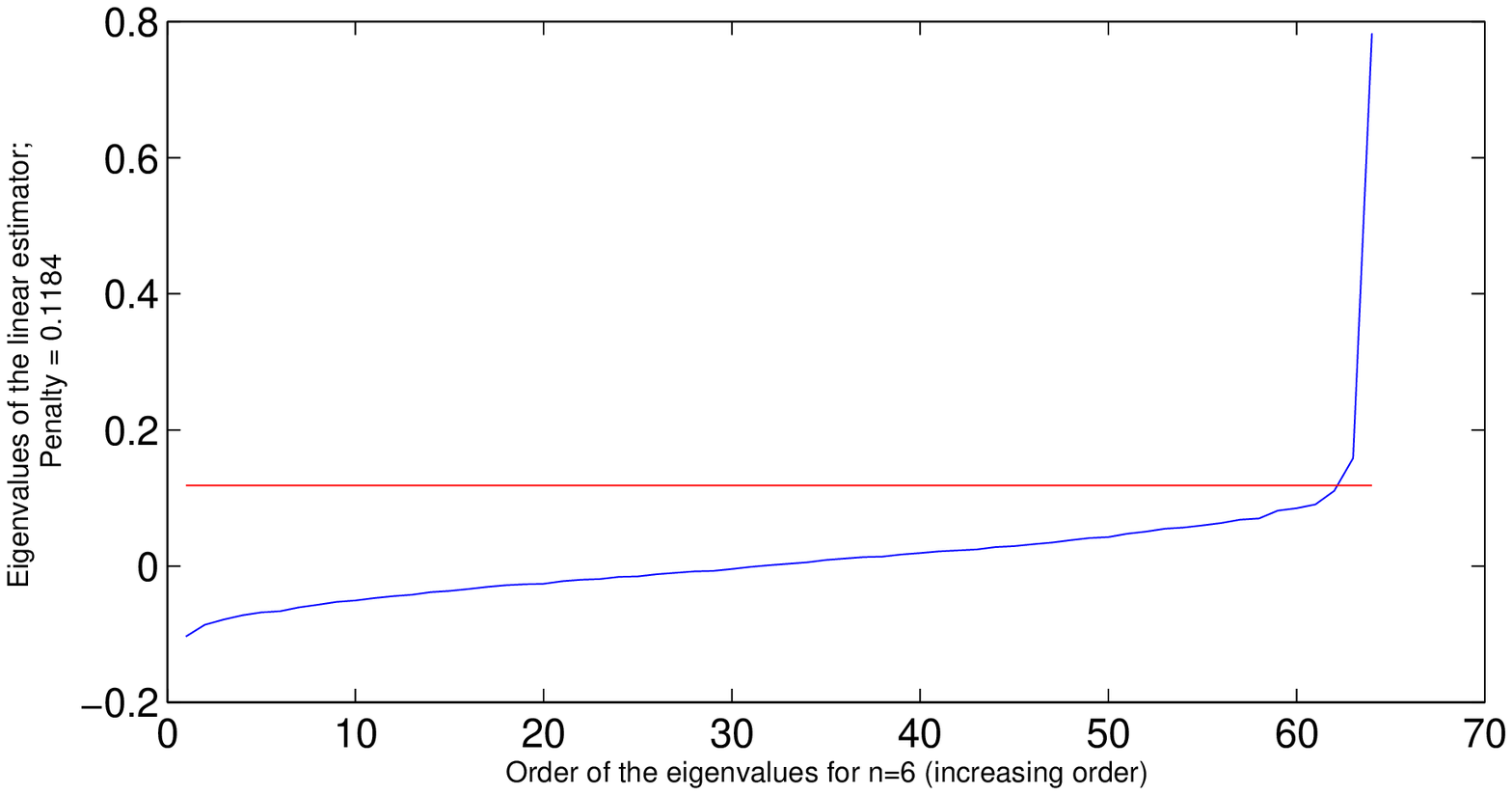}
\label{figeigenvalues}
\end{figure}

We have applied the method to real data sets concerning systems of 4 to 6 ions, which are Smolin states further manipulated. In Figure~\ref{figeigenvalues} we plot the eigenvalues of the linear estimator and the threshold given by the penalty. In each case, the method selects a rank equal to 2.

\newpage

\section{Conclusions}

We present here a method for reconstructing the quantum state of a system of $n$ qubits from all measurements, each repeated $m$ times. Such an experiment produce a huge amount of data to exploit in efficient way.

We revisit the inversion method and write an explicit formula for what is here called the linear estimator. This procedure does not produce a proper quantum state and has other well-known inconvenients. We consider projection of this state on the subspace of matrices with fixed rank and give an algorithm to select from data the rank which best suits the given quantum system. The method is very fast, as it comes down to choosing the eigenvalues larger than some threshold, which also appears in the penalty term. This threshold is of the same order as the error of the linear estimator. Its computation is crucial for good selection of the correct rank and it can be time consuming. Our algorithm also provides a consistent estimator of the true rank of the quantum system.

Our theoretical results provide a penalty term $\nu$ which has good asymptotic properties but our numerical results show that it is too large for most examples. Therefore we give an idea about how to evaluate closer the threshold by Monte-Carlo computation. This step can be time consuming but we can still improve on numerical efficiency (parallel computing, etc.).

In practice, the method works very well for large systems of small ranks, with significant eigenvalues. Indeed, there is a trade-off between the amount of data which will give small estimation error (and threshold) and the smallest eigenvalue that can be detected above this threshold. Neglecting eigenvalues comes down to reducing the number of parameters to estimate and reducing the variance, whereas large rank will increase the number of parameters and reduce the estimation bias.

{\bf Acknowledgements:} We are most grateful to M\u ad\u alin Gu\c t\u a and to Thomas Monz
for useful discussion and for providing us the experimental data used in this manuscript.


\section{Appendix}

{\bf Proof of Proposition~\ref{inversion}}
Actually, we can compute
\begin{eqnarray*}
(\mathbf{P}^{T} \mathbf{P})_{\b_1,\b_2} &=&
 \sum_{(\mathbf{r},\a)}
 \prod_{j\not\in E_{\b_1}} r_j \, \mathbf{I}{(a_j=b_{1,j})}
 \prod_{k\not\in E_{\b_2}} r_k \, \mathbf{I}{(a_k=b_{2,k})} .
\end{eqnarray*}
In case $\b_1=\b_2=\b$, we have
\begin{eqnarray*}
(\mathbf{P}^{T} \mathbf{P})_{\b,\b} &=&
\sum_{(\mathbf{r},\a)}
\left( \prod_{j\not\in E_{\b}} r_j \, \mathbf{I}{(a_j=b_j)} \right)^2
  \\
 & = & \sum_{(\mathbf{r},\a)}
 \prod_{j\not\in E_{\b}}  \mathbf{I}{(a_j=b_j)}
= 3^{d(\b)} 2^n.
\end{eqnarray*}
In case $\b_1 \not = \b_2$, we have either $E_{\b_1} = E_{\b_2}$ or $E_{\b_1}
\not = E_{\b_2}$.
If we suppose $E_{\b_1} = E_{\b_2}$,
$$
 \prod_{j\not\in E_{\b_1}} r_j \, \mathbf{I}{(a_j=b_{1,j})}
 \prod_{k\not\in E_{\b_2}} r_k \, \mathbf{I}{(a_k=b_{2,k})}
=0.
$$
Indeed, if this is not 0 it means $\a = \b_1=\b_2$ outside the set $E_{b_1}$,
that is $\b_1 = \b_2$ which contradicts our assumption.

If we suppose $E_{\b_1} \not = E_{\b_2}$, we have either $\b_1 \not = \b_2$ on
the set
$E_{\b_1}^C \cap E_{\b_1}^C$ and in this case one indicator in the product is
bound to be 0, or we have $\b_1 \not \not = \b_2$ on the set
$E_{\b_1}^C \cap E_{\b_1}^C$. In this last case, take $j_0$ in the symmetric
difference of sets
$E_{\b_1} \Delta E_{\b_2}$. Then,
\begin{eqnarray*}
&& (\mathbf{P}^{T} \mathbf{P})_{\b_1,\b_2} \\
&=&
 \sum_{(\mathbf{r},\a)}
 \prod_{j\not\in E_{\b_1}} r_j \, \mathbf{I}{(a_j=b_{1,j})}
 \prod_{k\not\in E_{\b_2}} r_k \, \mathbf{I}{(a_k=b_{2,k})} \\
&=& \sum_{(\mathbf{r},\a)}
 \prod_{j\not\in E_{\b_1}} \mathbf{I}{(a_j=b_{1,j})}
 \prod_{k\not\in E_{\b_2}} \mathbf{I}{(a_k=b_{2,k})} \prod_{j \in E_{\b_1}
\Delta E_{\b_2}} r_j \\
&=& \sum_{r_{j_0}\in \{-1,1\}} r_{j_0}
\sum_{\mathbf{r} \not = r_{j_0}}\sum_{ \a}
 \prod_{j\not\in E_{\b_1}} \mathbf{I}{(a_j=b_{1,j})} \\
&& \prod_{k\not\in E_{\b_2}} \mathbf{I}{(a_k=b_{2,k})} \prod_{j \in E_{\b_1}
\Delta E_{\b_2} \slash j_0} r_j = 0.
\end{eqnarray*}
\hfill $\Box$

{\bf Proof of Proposition \ref{lambda}}
It is easy to see that $\hat{\vec{\rho}}$ is an unbiased estimator.
We write its variance as follows:
\begin{eqnarray*}
&& Var(\hat \rho_{\b}) \\
&=&  \frac{1}{3^{2 d(\b)} 4^n}
\sum_{\mathbf{a} \in  \mathcal{E}^{n} } Var \left( \sum_{\r \in \mathcal{R}^n}
\frac 1m \sum_{i=1}^m \delta_{R^{\mathbf{a},i}, \mathbf{r}}
\mathbf{P}_{(\mathbf{r},\mathbf{a}),\b} \right)\\
&=&  \frac{1}{3^{2 d(\b)} 4^n m^2}
\sum_{\mathbf{a} \in  \mathcal{E}^{n} }
 \sum_{\r \in \mathcal{R}^n} m p_{(\r, \a)}
\mathbf{P}_{(\mathbf{r},\mathbf{a}),\b}^2 \\
&&  - \frac{1}{3^{2 d(\b)} 4^n m^2}
\sum_{\mathbf{a} \in  \mathcal{E}^{n} } m \left( \sum_{\r \in \mathcal{R}^n}  p_{(\r, \a)}
\mathbf{P}_{(\mathbf{r},\mathbf{a}),\b}\right)^2
\\
& = & \frac{1}{3^{2 d(\b)} 4^n m}
\sum_{(\mathbf{r},\mathbf{a})\in(\mathcal{R}^n\times \mathcal{E}^{n})}
 p_{(\r, \a)}\prod_{j\not\in E_{\b}}  \mathbf{I}{(a_j=b_j)}\\
&& -\frac 1m \sum_{\mathbf{a} \in  \mathcal{E}^{n} }\left( \frac{1}{3^{d(\b)}
2^n}
 \sum_{\r \in \mathcal{R}^n} p_{(\r, \a)} \prod_{j\not\in E_{\b}} r_j
\mathbf{I}{(a_j=b_j)}\right)^2\\
&\leq & \frac{1}{3^{ d(\b)} 4^n m} .
\end{eqnarray*}

Finally, let us prove the last point.
We will use the following result due to \cite{Tropp}.
\begin{thm}[Matrix Hoeffding's inequality \cite{Tropp}]
\label{thmTropp}
 Let $X_1$, ..., $X_p$ be independent centered self-adjoint random
 matrices with values in $\mathbb{C}^{d\times d}$, and let us assume
that there are deterministic self-adjoint matrices $A_1$, ..., $A_p$
such that, for all $i_in\{1,...,p\}$, $A_i^2-X_i^2$ is a.s. nonnegative.
Then, for all $t>0$,
$$
\mathbb{P}\left( \biggl\|\sum_{i=1}^{p} X_i \biggr\|^{2} \geq t \right)
 \leq d\exp\left(\frac{-t^2}{8\sigma^2}\right)
$$
where $\sigma^{2}=\|\sum_{k=1}^{p} A_k^2\|$.
\end{thm}
We have:
\begin{align*}
 \hat{\rho} - \rho
   & = \sum_\mathbf{b} (\hat{\rho}_\mathbf{b}-\rho_\mathbf{b})\sigma_{\mathbf{b}} \\
   & = \sum_\mathbf{b} \sum_\mathbf{r} \sum_\mathbf{a} \frac{\mathbf{P}_{(\mathbf{r},\mathbf{a}),\mathbf{b}}}
                            {3^{d(\mathbf{b})}2^n}(\hat{p}_{\mathbf{r},\mathbf{a}}
                  -p_{\mathbf{r},\mathbf{a}})\sigma_{\mathbf{b}} \\
   & = \sum_\mathbf{b} \sum_\mathbf{r} \sum_\mathbf{a} \sum_{i} \frac{\mathbf{P}_{(\mathbf{r},\mathbf{a}),\mathbf{b}}}
                            {3^{d(\mathbf{b})}2^n m}(\mathds{1}_{R^{i,\mathbf{a}}=\mathbf{r}}
                  -p_{\mathbf{r},\mathbf{a}})\sigma_{\mathbf{b}} \\
   & = \sum_\mathbf{a} \sum_{i} \underbrace{\sum_\mathbf{b} \sum_\mathbf{r} \frac{\mathbf{P}_{(\mathbf{r},\mathbf{a}),\mathbf{b}}}
                            {3^{d(\mathbf{b})}2^n}(\mathds{1}_{R^{i,\mathbf{a}}=\mathbf{r}}
                  -p_{\mathbf{r},\mathbf{a}})\sigma_{\mathbf{b}}}_{=:X_{i,\mathbf{a}}}.
\end{align*}
Note that the $X_{i,\mathbf{a}}$, for $i\in\{1,...,m\}$ and $\mathbf{a}\in\mathcal{E}^{n}$, are iid
self-adjoint centered random matrices. Moreover, we have:
\begin{align*}
 \| X_{i,\mathbf{a}} \| & = \left\| \sum_\mathbf{b} \sum_\mathbf{r} \frac{\mathbf{P}_{(\mathbf{r},\mathbf{a}),\mathbf{b}}}
                            {3^{d(\mathbf{b})}2^n m}(\mathds{1}_{R^{i,\mathbf{a}}=\mathbf{r}}
                  -p_{\mathbf{r},\mathbf{a}})\sigma_{\mathbf{b}} \right\| \\
  & \leq  \sum_\mathbf{b} \sum_\mathbf{r} \left| \frac{\mathbf{P}_{(\mathbf{r},\mathbf{a}),\mathbf{b}}}
                            {3^{d(\mathbf{b})}2^n m}\right| \left| \mathds{1}_{R^{i,\mathbf{a}}=\mathbf{r}}
                  -p_{\mathbf{r},\mathbf{a}} \right| \underbrace{\left\| \sigma_{\mathbf{b}} \right\|}_{=1} \\
  &  = \sum_\mathbf{b}  \left| \frac{\mathbf{P}_{(\mathbf{r},\mathbf{a}),\mathbf{b}}}
                            {3^{d(\mathbf{b})}2^n m}\right| \underbrace{\sum_\mathbf{r} \left| \mathds{1}_{R^{i,\mathbf{a}}=\mathbf{r}}
                  -p_{\mathbf{r},\mathbf{a}} \right|}_{\leq 2} \\
  &  = \frac{2}{2^n m} \sum_\mathbf{b} \frac{1}{3^{d(\mathbf{b})}} \prod_{j\notin E_{\mathbf{b}}} \mathds{1}_{a_j = b_j} \\
  &  \leq \frac{2}{2^n m} \sum_{\ell=0}^{n} \sum_{
\tiny{
\begin{array}{c}
\mathbf{b} \text{ such that} \\
 d(\mathbf{b})=\ell \\
 \forall j\notin E_{\mathbf{b}}, a_j = b_j
\end{array}
}
}
 \frac{1}{3^{\ell}} \\
  & = \frac{2}{2^n m} \sum_{\ell=0}^{n} {\ell \choose n} \frac{1}{3^{\ell}}
= \frac{2}{2^n m} \left(1+\frac{1}{3}\right)^{n}  = \frac{2}{m} \left(\frac{2}{3}\right)^n.
\end{align*}
This proves that $A_{i,\mathbf{a}}^2 - X_{i,\mathbf{a}}^2$ is nonnegative where
$A_{i,\mathbf{a}} = \frac{2}{m}\left(\frac{2}{3}\right)^n I$. So we can apply
Theorem~\ref{thmTropp}, we have:
$$\sigma^{2}=\|\sum_{i,\mathbf{a}} A_{i,\mathbf{a}}^2\| = \frac{4}{m}\left(\frac{4}{3}\right)^n $$
and so
\begin{eqnarray*}
&&\mathbb{P}\left( \left\| \hat{\rho} - \rho \right\|^{2} \geq t \right)
 = \mathbb{P}\left( \biggl\|\sum_{i,\mathbf{a}} X_{i,\mathbf{a}} \biggr\|^{2} \geq t \right)\\
& \leq &2^n \exp\left(\frac{-t^2 m}{32} \left(\frac{3}{4}\right)^n \right).
\end{eqnarray*}
We put
$$ \varepsilon = 2^n \exp\left(\frac{-t^2 m}{32} \left(\frac{3}{4}\right)^n \right), $$
this leads to:
$$
\mathbb{P}\left( \left\| \hat{\rho} - \rho \right\|^{2} \geq
4\sqrt{2 \left(\frac{4}{3}\right)^n   \frac{n\log(2) - \log(\varepsilon) }{m} }\right) \leq \varepsilon.
$$
\hfill $\Box$

{\bf Proof of Theorem~\ref{res}} From the definition (\ref{pen}) of our
estimator, we have, for any Hermitian, positive semi-definite matrix $R$,
\begin{eqnarray*}
\left\|\hat{\rho}_{\nu} - \hat{\rho} \right\|_F^{2} + \nu
{\rm rank}(\hat{\rho}_{\nu})
& \leq & \left\|R - \hat{\rho}  \right\|_F^{2} + \nu {\rm rank}(R).
\end{eqnarray*}
We deduce that
\begin{eqnarray*}
&&\left\|\hat{\rho}_{\nu} - {\rho} \right\|_F^{2}\\
& \leq & \left\| R - {\rho} \right\|_F^{2}
+ 2 {\rm Tr}((\hat{\rho} - \rho)^\star(R -
\hat{\rho}_{\nu}))\\
&& + \nu ({\rm rank}(R) - {\rm rank}(\hat{\rho}_{\nu})) \\
&\leq & \left\| R - {\rho} \right\|_F^{2} + 2 \nu {\rm rank}(R)
+2 \|\hat{\rho} - \rho\| \times \|R -
\hat{\rho}_{\nu}\|_1\\
&& - \nu ({\rm rank}(R) + {\rm rank}(\hat{\rho}_{\nu})).
\end{eqnarray*}
Further on, we have
\begin{eqnarray*}
&& \|R - \hat{\rho}_{\nu}\|_1 \\
&\leq & ({\rm rank}(R)+ {\rm rank}(\hat{\rho}_{\nu}))^{1/2} \|R -
\hat{\rho}_{\nu}\|_F\\
&\leq & ({\rm rank}(R)+ {\rm rank}(\hat{\rho}_{\nu}))^{1/2}
(\|\rho - \hat{\rho}_{\nu}\|_F+\|R - \rho\|_F)
\end{eqnarray*}
We apply two times the inequality $2 A\cdot B \leq \epsilon A^2+
\epsilon^{-1}B^2$ for any real numbers
$A,\, B$ and $\epsilon >0$. We actually use $\epsilon = 1+ \theta /2$ and
$\epsilon = \theta /2$, respectively, and get
\begin{eqnarray*}
&& \left\|\hat{\rho}_{\nu} - {\rho} \right\|_F^{2} \\
&\leq & \left\| R - {\rho} \right\|_F^{2} + 2 \nu {\rm rank}(R)
- \nu ({\rm rank}(R) + {\rm rank}(\hat{\rho}_{\nu}))\\
&& +(1+\theta)({\rm rank}(R)+ {\rm
rank}(\hat{\rho}_{\nu}))\|\hat{\rho} - \rho\|^2\\
&& + (1+\frac{\theta }2 )^{-1} \left\|\hat{\rho}_{\nu} - {\rho}
\right\|_F^{2} + (\frac \theta 2)^{-1}\|R - \rho\|_F^2.
\end{eqnarray*}
By rearranging the previous terms, we get that for any Hermitian matrix $R$
$$
\left\|\hat{\rho}_{\nu} - {\rho} \right\|_F^{2}
\leq c^2(\theta) \|R - \rho\|_F^2 + 2 c(\theta)  \nu {\rm rank}(R),
$$
provided that $\nu \geq (1+\theta) \|\hat{\rho} - \rho\|^2$.
By following \cite{BSW}, the least possible value for $\|R - \rho\|_F^2$
is $\sum_{j>k} \lambda_j^2(\rho)$ if the matrices $R$ have rank $k$.
Moreover, this value is obviously attained by the projection of $\rho$ on
the space of the eigenvectors associated to the $k$ largest eigenvalues. This
helps us conclude the proof of the theorem.
\hfill $\Box$

{\bf Proof of Corollary~\ref{hatrank}} Recall that $\hat k$ is the
largest $k$ such that $\lambda_k(\hat{\rho}) \geq \sqrt{\nu}$. We have
$$
\mathbb{P}(\hat k \not = k)  = \mathbb{P}(\lambda_k(\hat{\rho})\leq
\sqrt{\nu} \mbox{ or } \lambda_{k+1}(\hat{\rho}) \geq \sqrt{\nu}).
$$
Now, $\lambda_k({\rho}) \leq \lambda_k(\hat{\rho})+ \|\hat{\rho}
- \rho\|$ and $\lambda_{k+1}({\rho})\geq \lambda_{k+1}(\hat{
\rho})-\|\hat{\rho} - \rho\|$. Thus,
$$
\mathbb{P}(\hat k \not = k) \leq \mathbb{P}(\|\hat{\rho} - \rho\|\geq
\min\{\lambda_k({\rho})-\sqrt{\nu}, \sqrt{\nu}-\lambda_{k+1}({
\rho})\})
$$
and this is smaller than $\mathbb{P}(\|\hat{\rho} - \rho\|\geq \delta
\sqrt{\nu})$, by the assumptions of the Corollary.
\hfill $\Box$


\begin{thebibliography}{0}
\expandafter\ifx\csname natexlab\endcsname\relax\def\natexlab#1{#1}\fi
\expandafter\ifx\csname bibnamefont\endcsname\relax
  \def\bibnamefont#1{#1}\fi
\expandafter\ifx\csname bibfnamefont\endcsname\relax
  \def\bibfnamefont#1{#1}\fi
\expandafter\ifx\csname citenamefont\endcsname\relax
  \def\citenamefont#1{#1}\fi
\expandafter\ifx\csname url\endcsname\relax
  \def\url#1{\texttt{#1}}\fi
\expandafter\ifx\csname urlprefix\endcsname\relax\def\urlprefix{URL }\fi
\providecommand{\bibinfo}[2]{#2}
\providecommand{\eprint}[2][]{\url{#2}}

\end{thebibliography}


\begin{thebibliography}{99}

\bibitem{AF+}
{ M. Asorey, P. Facchi, G. Florio, V. I. Man'ko, G. Marmo, S. Pascazio and
 E. C. G. Sudarshan}
Phys. Lett. A {\bf 375}, 861 (2011).

\bibitem{AS}
{K. M. R. Audenaert and S. Scheel},
New J. Phys. {\bf 11}, 023028 (2009).

\bibitem{BBG}
{E. Bagan, M. A. Ballester, R. D. Gill, A. Monras, and R.
Mu\~noz Tapia}
Phys. Rev. A {\bf 73}, 032301 (2006).

\bibitem{BAP}
{K. Banaszek, G. M. D'Ariano, M. G. A. Paris, and M. F.
Sacchi} (1999)
Phys. Rev. A {\bf 61}, 010304 (1999).

\bibitem{BSG}
{J. T. Barreiro, P. Schindler, O. Guhne, T. Monz,
M. Chwalla, C. F. Roos, M. Hennrich, and R. Blatt},
Nature Phys. {\bf 6}, 943 (2010).

\bibitem{BK1}
{R. Blume-Kohout}
New J. Phys. {\bf 12}, 043034 (2010).

\bibitem{BK2}
{R. Blume-Kohout}
Hedged Maximum Likelihood Quantum State Estimation,
Phys. Rev. Lett. {\bf 105}, 200504 (2010).

\bibitem{BD+}
{G. Brida, I. P. Degiovanni, A. Florio, M. Genovese, P. Giorda, A. Meda,
M. G. A. Paris and A. Shurupov}
Phys. Rev. Lett. {\bf 104}, 100501 (2010).

\bibitem{BSW}
{F. Bunea, Y. She, and M. H. Wegkamp} (2011)
Ann. Statist., {\bf 39}(2), pp. 1282-1309.

\bibitem{Candes}
{E. J. Cand\`es and Y. Plan} (2010)
arxiv:1101.0339 (math.ST)

\bibitem{Filippov}
{S. N. Filippov and V. I. Man'ko}
Physica Scripta {\bf T143}, 014010 (2011).

\bibitem{GLF}
{D. Gross, Y.-K. Liu, S. T. Flammia, S. Becker and J. Eisert}
{Phys. Rev. Lett.},{\bf 105}, 150401 (2010).

\bibitem{Gross}
{D. Gross}
{IEEE Trans. on Information Theory}, {\bf 57}, 1548-1566  (2011).

\bibitem{GKD}
{M. Gu\c t\u a, T. Kypraios and I. Dryden} (2012)
{New J. Physics}, {\bf 14}, 105002.

\bibitem{Ivonovic}
{I. D. Ivonovic}
J. Phys. A {\bf 14}, 3241 (1981).

\bibitem{JKM}
{D. F. V. James, P. G. Kwiat, W. J. Munro, and A. G. White}
Phys. Rev. A {\bf 64}, 052312  (2001).

\bibitem{Klopp}
{O. Klopp}
{Electronic J. Statist.}, {\bf 5}, pp. 1161-1183  (2011).

\bibitem{Kolt}
{V. Koltchinskii}
{Ann. Statist.}, {\bf 39}(6), pp. 2936-2973  (2011).

\bibitem{KLT}
{V. Koltchinskii, K. Lounici and A. B. Tsybakov}
{Ann. Statist.}, {\bf 39}(5), pp. 2302-2329 (2011).

\bibitem{KoltStF}
{V. Koltchinskii}
{\it Oracle Inequalities in Empirical Risk Minimization and Sparse Recovery Problems},
Ecole d'\'et\'e de Probabilit\'es de Saint-Flour XXXVIII, Springer Lecture Notes in
Mathematics  (2011).

\bibitem{L}
{A. I. Lvovsky}
J. Opt. B. {\bf 6}(6), S556 (2004).

\bibitem{Monz14}
{T. Monz, Ph. Schindler, J. Barreiro, M. Chwalla, D. Nigg, W. Coish,
M. Harlander, W. H\"ansel, M. Hennrich, and R. Blatt}.
arXiv:1009.6126v2 (quant-ph)  (2011)

\bibitem{NW}
{S. Negahban and M. J. Wainwright}
arxiv:0912.5100 (math.ST), to appear in {\it Ann. Statist.}  (2009)

\bibitem{RMH}
{J. \v{R}eh\'a\v{c}ek, D. Mogilevtsev and Z. Hradil}
Phys. Rev. Lett. {\bf 105}, 010402 (2010).

\bibitem{RV98}
{G. C. Reinsel and R. P. Velu}.
Multivariate Reduced-Rank Regression: Theory and Applications.
{\it Lecture Notes in Statistics, Springer, New York} (1998).

\bibitem{RT}
{A. Rohde and A. B. Tsybakov}
{Ann. Statist.}, {\bf 39}(2), pp. 887-930 (2011).

\bibitem{Tropp}
{J. A. Tropp}
{Foundations of Comput. Mathem.}  (2011).

\bibitem{VR}
{K. Vogel and H. Risken}
Phys. Rev. A {\bf 40}, 2847
(1989).

\bibitem{Y}
{K. Yamagata}
Int. J. Quantum Inform. {\bf 9}(4), 1167 (2011).

\end{thebibliography}
\end{document}